\pdfoutput=1
\documentclass[11pt]{article}
\usepackage[left=1in,right=1in,top=1in,bottom=1in]{geometry}
\usepackage{times}
\usepackage{expl3}
\usepackage{cite}
\usepackage[table]{xcolor}
\usepackage{multirow}
\usepackage{stackengine} 
\usepackage{hhline}
\usepackage{lipsum}
\usepackage{titlesec}
\usepackage{wrapfig}
\usepackage{enumerate}
\usepackage{epsfig}
\usepackage{amsmath}
\usepackage{tabularx}
\usepackage{array}
\usepackage{booktabs}
\usepackage{enumitem}
\usepackage{bbm}
\usepackage{calc}
\usepackage{graphicx}
\usepackage{amsmath}
\usepackage[title]{appendix}
\usepackage{amssymb}
\usepackage{epstopdf}
\usepackage{boldline}
\usepackage{arydshln}
\usepackage{calligra}
\usepackage{bm}
\usepackage{url}
\usepackage{blindtext}
\usepackage{accents}
\usepackage{amsthm}

\newcommand{\define}{\stackrel{\mbox{\tiny def}}{=}}

\newtheorem{definition}{Definition}
\newtheorem{theorem}{Theorem}

\newtheorem{corollary}{Corollary}
\newtheorem{lemma}{Lemma}
\newtheorem{example}{Example}

\usepackage{mathtools}
\usepackage{epstopdf}
\usepackage{balance}
\usepackage{thmtools}
\usepackage{thm-restate}
\usepackage{hyperref}
\usepackage{cleveref}
\usepackage[mathscr]{euscript}

\usepackage[ruled,vlined]{algorithm2e}
\include{pythonlisting}

\newcommand{\ostar}{\mathbin{\mathpalette\make@circled\star}}

\makeatletter
\newcommand{\removelatexerror}{\let\@latex@error\@gobble}
\makeatother
\makeatletter
\newcommand*{\rom}[1]{\expandafter\@slowromancap\romannumeral #1@}
\makeatother

\ExplSyntaxOn
\newcommand\latinabbrev[1]{
  \peek_meaning:NTF . {
    #1\@}%
  { \peek_catcode:NTF a {
      #1.\@ }%
    {#1.\@}}}
\ExplSyntaxOff

\titleclass{\subsubsubsection}{straight}[\subsubsection]

\begin{document}
\vspace{1cm}
\title{Spectral and Nilpotent Matrix Orderings: Comparison and Applications in Dynamic Systems}
\vspace{1.8cm}
\author{Shih-Yu~Chang
\thanks{Shih-Yu Chang is with the Department of Applied Data Science,
San Jose State University, San Jose, CA, U. S. A. (e-mail: {\tt
shihyu.chang@sjsu.edu}). 
           }}

\maketitle

\begin{abstract}
In our earlier work, we proposed the \emph{Spectral and Nilpotent Ordering} (SNO) as a new framework that extends matrix comparison beyond the Hermitian setting by incorporating both spectral and nilpotent structures. Building on that foundation, the present paper develops concrete certificates and applications of SNO. First, we employ generalized Gershgorin theorems to design certificates for spectral ordering that avoid direct eigenvalue computation and analyze their robustness under perturbations. Second, we introduce rank-based criteria that provide certificates for ordering the nilpotent parts of matrices without requiring a full Jordan decomposition. Finally, we apply the SNO framework to linear dynamical systems, where we construct a hierarchy of stability orderings that capture both asymptotic and transient behaviors. These contributions advance the theoretical underpinnings of SNO and demonstrate its potential as a versatile tool for operator analysis, computational methods, and stability studies in complex systems.
\end{abstract}
\begin{keywords}
Spectral and Nilpotent Ordering (SNO); Matrix comparison; Spectral ordering; Nilpotent ordering; Gershgorin theorem; Stability analysis; Linear dynamical systems.
\end{keywords}

\section{Introduction}\label{sec: Introduction}

Matrices are indispensable tools across pure mathematics and engineering disciplines. They model linear transformations, encode system dynamics, and underpin numerical methods such as decompositions (e.g., LU, QR, generalized inverses), stability analysis, and perturbation theory. In engineering, matrices represent state-space systems, control designs, signal processing operations, and quantum feedback circuits. Recent advances in the field have strengthened the theoretical foundations of matrix decompositions and perturbation bounds \cite{WeiGLUC2025,xu2025lu}, while other studies have expanded applications to the control of quantum linear systems and the analysis of networked dynamical structures \cite{ZhangQuantum2015,ZhangMMatrices2010}.

Matrix comparison plays a central role in mathematics, science, and engineering, as it provides a rigorous framework for analyzing similarities and differences between linear operators. In mathematics, spectral and structural comparisons are essential for understanding stability, perturbation theory, and functional calculus. In scientific computing, efficient matrix comparison methods enable dimensionality reduction, pattern recognition, and tensor-based analysis of high-dimensional data. In engineering applications, particularly in control theory and signal processing, comparing matrices under spectral or nilpotent criteria helps assess system stability, robustness, and performance under uncertainty. Such techniques highlight the interdisciplinary importance of matrix comparison, bridging theoretical insights with practical problem solving~\cite{bhatia1997matrix}.

The L\"owner ordering provides a partial order structure to compare matrices with the same size in the space of Hermitian matrices (or more generally, self-adjoint operators). If we are given two Hermitian matrices $\bm{A}$ and $\bm{B}$, we express $\bm{A} \preceq \bm{B}$ in the sense of L\"owner ordering if and only if $\bm{B} - \bm{A}$ is positive semidefinite. This ordering plays an important role in operator theory and matrix analysis, as it offers a rigorous framework for comparing spectra and establishing monotonicity of operator functions. In applications, the L\"owner ordering is widely used in quantum information theory, optimization, and statistics, where positivity and dominance relations between covariance matrices or density operators are critical~\cite{ando1979concavity}. Several works about hypercomplex analysis are based on L\"owner ordering~\cite{chang2024generalizedCDJ,chang2024generalizedJensen,chang2024multivariate}.

Classical L\"owner ordering provides a powerful framework for comparing Hermitian matrices by exploiting their real spectra, but it inherently excludes non-Hermitian systems, whose eigenvalues may be complex. In our recent work~\cite{chang2025matrixSNO}, we extend the spirit of L\"owner ordering to a broader setting by introducing the \emph{Spectral and Nilpotent Ordering} (SNO). Unlike the L\"owner ordering, which is restricted to Hermitian operators, the SNO framework defines a total ordering relation for complex numbers, enabling the comparison of matrices with complex spectra. Moreover, by incorporating dominance relations for the nilpotent parts of Jordan blocks, SNO systematically addresses the structure of non-diagonalizable matrices. In this way, our approach generalizes L\"owner ordering beyond the Hermitian domain, providing novel monotonicity and convexity results for matrix functions and significantly broadening the scope of matrix comparison theory.

This paper advances the theory of Spectral and Nilpotent Ordering (SNO) by addressing several open challenges in matrix comparison. A key difficulty lies in comparing two matrices without fully computing their Jordan decompositions. To this end, we design certificates for spectral ordering based on a generalized Gershgorin theorem, which allow us to determine ordering relations without directly computing eigenvalues. We further investigate the robustness of spectral ordering under perturbations, providing tools for practical applications where exact spectra are inaccessible. On the nilpotent side, we establish certificates using matrix rank relations to characterize the ordering of nilpotent parts when two matrices share identical spectra, thereby avoiding the explicit computation of canonical forms. Finally, we demonstrate how SNO can be applied to stability analysis in linear dynamical systems. We introduce a hierarchy of stability orderings within this framework and prove several theorems that reveal how SNO captures both asymptotic and transient stability behaviors, offering a new lens for studying system dynamics.

The paper is organized as follows. Section~\ref{sec:Review Spectral and Nilpotent Ordering} provides a review of spectral and nilpotent ordering. In Section~\ref{sec:Spectral Ordering}, we present certificates for spectral ordering and discuss their role in comparing matrices through eigenvalue structures. Section~\ref{sec:Nilpotent Part Ordering} develops certificates for nilpotent ordering, emphasizing the influence of Jordan block structures on ordering relations. Section~\ref{sec:SNO Applications} explores applications of spectral and nilpotent ordering to the stability analysis of dynamical systems, highlighting the utility of the SNO framework in assessing comparative stability. 

\section{Review Spectral and Nilpotent Ordering}\label{sec:Review Spectral and Nilpotent Ordering}

In this section, we will provide a quick reviw about Spectral and Nilpotent Ordering (SNO) given by~\cite{chang2025matrixSNO}. We begin by defining \emph{dominance order} between two vectors of natural numbers. Let \( m,n \) be two natural numbers, and let \( \bm{p}= (p_1, p_2, \ldots, p_k) \) and \( \bm{q} = (q_1, q_2, \ldots, q_l) \) be two partitions of $m$ and $n$, respectively, written in non-increasing order:
\[
p_1 \geq p_2 \geq \cdots \geq p_k > 0, \quad q_1 \geq q_2 \geq \cdots \geq q_l > 0.
\]

We say that $\bm{p}$ precedes $\bm{q}$ in the dominance order, written $\bm{p}\trianglelefteq\bm{q}$, if the following conditions are satisfied:\\

1. Partial sums condition:
\begin{eqnarray}\label{eq1: dom ordering}
   \sum_{i=1}^j p_i \leq \sum_{i=1}^j q_i \quad \text{for all } j \geq 1,
\end{eqnarray}
where parts beyond the length of a partition are treated as \(0\) (e.g., if \(\bm{p}\) has 3 parts but \(\bm{q}\) has 4 parts, then \(p_4 = 0\)). \\

2. Total sum condition:
   \[
   \sum_{i=1}^k p_i = m, \quad \sum_{i=1}^l q_i = n,
   \]
   ensuring both are partitions of numbers $m$ and $n$.  Note that $\bm{p}\triangleleft\bm{q}$ indicates that  $ \sum_{i=1}^j p_i < \sum_{i=1}^j q_i \quad \text{for all } j \geq 1$.

Given two matrices $\bm{X}_1$ and $\bm{X}_2$ with the same dimensions $m \times m$ and Jordan decomposition, these two matrices $\bm{X}_1$ and $\bm{X}_2$ can be expressed by
\begin{eqnarray}\label{eq: matrices setup}
\bm{X}_1&=& \bm{U}_1\left(\bigoplus\limits_{k_1=1}^{K_1}\bigoplus\limits_{i_1=1}^{\alpha_{k_1}^{(\mathrm{G})}}\bm{J}_{m_{k_1,i_1}}(\lambda_{k_1})\right)\bm{U}_1^{-1}, \nonumber \\
\bm{X}_2&=& \bm{U}_2\left(\bigoplus\limits_{k_2=1}^{K_2}\bigoplus\limits_{i_2=1}^{\alpha_{k_2}^{(\mathrm{G})}}\bm{J}_{m_{k_2,i_2}}(\lambda_{k_2})\right)\bm{U}_2^{-1},
\end{eqnarray}   
where $\alpha_{k_1}^{(\mathrm{G})}$ is the geometry multiplicities for the eigenvalue $\lambda_{k_1}$ of the matrix $\bm{X}_1$ and $\alpha_{k_2}^{(\mathrm{G})}$ is the geometry multiplicities for the eigenvalue $\lambda_{k_2}$ of the matrix $\bm{X}_2$, respectively.

We have the following definition to compare two matrices $\bm{X}_1$ and $\bm{X}_2$ by their representations $\mathfrak{R}(\bm{X}_1)$ and $\mathfrak{R}(\bm{X}_2)$. Let $\alpha_k^{(\mathrm{A})}(\bm{X})$ is the algebraic multiplicity of the eigenvalue $\lambda_K$ with respect to the matrix $\bm{X}$. 
\begin{definition}\label{def: SNO}
If we set $[\underbrace{\lambda_k(\bm{X}),\ldots,\lambda_k(\bm{X})}_{\alpha_k^{(\mathrm{A})}(\bm{X})}]\define\bm{\lambda_k}(\bm{X})$ and $[\underbrace{m_{k,1}(\bm{X}),\ldots,m_{k,\alpha_{k}^{(\mathrm{G})}(\bm{X})}(\bm{X})}_{\sum\limits_{i=1}^{\alpha_k^{(\mathrm{G})}(\bm{X})}m_{k,i}(\bm{X})=\alpha_k^{(\mathrm{A})}(\bm{X})}]\define\bm{m}_k(\bm{X})$ for $k=1,2,\ldots,K$, we can represent matrices $\bm{X}_1$ and $\bm{X}_2$ as 
\begin{eqnarray}\label{eq1: def: SNO}
\mathfrak{R}(\bm{X}_1)&=&[\bm{\lambda_1}(\bm{X}_1),\ldots,\bm{\lambda}_{K_1}(\bm{X}_1),\bm{m}_1(\bm{X}_1),\ldots,\bm{m}_{K_1}(\bm{X}_1)], \nonumber \\
\mathfrak{R}(\bm{X}_2)&=&[\bm{\lambda_1}(\bm{X}_2),\ldots,\bm{\lambda}_{K_2}(\bm{X}_2),\bm{m}_1(\bm{X}_2),\ldots,\bm{m}_{K_2}(\bm{X}_2)],
\end{eqnarray}   
where we assume that both $\lambda_j(\bm{X}_1)$ (for $j \in 1,2,\ldots,K_1$) and $\lambda_j(\bm{X}_2)$ (for $j \in 1,2,\ldots,K_2$) are arranged in the descendig order in the ordering  sense of $\leq_{\text{lex}}$~\footnote{Two complex numbers $a\leq_{\text{lex}} b$ if $\Re(a) < \Re(b)$ or $\Re(a) = \Re(b)$ with $\Im(a) < \Im(b)$, where $\Re$ and $\Im$ are used to represent the real and the imaginary parts, respectively. }.

We say that $[\bm{m}_1(\bm{X}_1),\ldots,\bm{m}_{K_1}(\bm{X}_1)] \preceq_{\mbox{\tiny N}} [\bm{m}_1(\bm{X}_2),\ldots,\bm{m}_{K_2}(\bm{X}_2)]$, where $\preceq_{\mbox{\tiny N}}$ represents the order between $\max(K_1, K_2)$ distinct eigenvalues by considering their nilpotent structures~\footnote{If $K_1 \neq K_2$, the shorter list will be added $|K_1 - K_2|$ zero entries before applying $\trianglelefteq$ comparison.}, if and only if 
there is some $k \in \{1,2,\ldots, \max(K_1, K_2)\}$ such that $\bm{m}_j(\bm{X}_1)=\bm{m}_j(\bm{X}_2)$ for all $j < k$ and $\bm{m}_k(\bm{X}_1)\trianglelefteq\bm{m}_k(\bm{X}_2)$.  

We say that $\bm{X}_1 \preceq_{\mbox{\tiny SN}} \bm{X}_2$, where $\preceq_{\mbox{\tiny SN}}$ represents the order between matrices by considering their spectral and nilpotent structures, if and only if 
\begin{eqnarray}\label{eq3: def: SNO}
[\bm{\lambda_1}(\bm{X}_1),\ldots,\bm{\lambda}_{K_1}(\bm{X}_1)] \preceq_w  
[\bm{\lambda_1}(\bm{X}_2),\ldots,\bm{\lambda}_{K_2}(\bm{X}_2)] 
\end{eqnarray}   
or
\begin{eqnarray}\label{eq4: def: SNO}
[\bm{\lambda_1}(\bm{X}_1),\ldots,\bm{\lambda}_{K_1}(\bm{X}_1)] &=&  
[\bm{\lambda_1}(\bm{X}_2),\ldots,\bm{\lambda}_{K_2}(\bm{X}_2)] \nonumber \\
\mbox{with} & & [\bm{m}_1(\bm{X}_1),\ldots,\bm{m}_{K_1}(\bm{X}_1)] \preceq_{\mbox{\tiny N}} [\bm{m}_1(\bm{X}_2),\ldots,\bm{m}_{K_2}(\bm{X}_2)],
\end{eqnarray}   
where $\preceq_w$ represents the weak majorization between two vectors of complex numbers under $\leq_{\text{lex}}$.

Moreover, 
We say that $\bm{X}_1 \prec_{\mbox{\tiny SN}} \bm{X}_2$, where $\prec_{\mbox{\tiny SN}}$ represents the order between matrices by considering their spectral and nilpotent structures, if and only if 
\begin{eqnarray}\label{eq5: def: SNO}
[\bm{\lambda_1}(\bm{X}_1),\ldots,\bm{\lambda}_{K_1}(\bm{X}_1)] \prec_w  
[\bm{\lambda_1}(\bm{X}_2),\ldots,\bm{\lambda}_{K_2}(\bm{X}_2)] 
\end{eqnarray}   
or
\begin{eqnarray}\label{eq6: def: SNO}
[\bm{\lambda_1}(\bm{X}_1),\ldots,\bm{\lambda}_{K_1}(\bm{X}_1)] &=&  
[\bm{\lambda_1}(\bm{X}_2),\ldots,\bm{\lambda}_{K_2}(\bm{X}_2)] \nonumber \\
\mbox{with} & & [\bm{m}_1(\bm{X}_1),\ldots,\bm{m}_{K_1}(\bm{X}_1)] \prec_{\mbox{\tiny N}} [\bm{m}_1(\bm{X}_2),\ldots,\bm{m}_{K_2}(\bm{X}_2)].
\end{eqnarray}   
We say that $[\bm{m}_1(\bm{X}_1),\ldots,\bm{m}_{K_1}(\bm{X}_1)] \prec_{\mbox{\tiny N}} [\bm{m}_1(\bm{X}_2),\ldots,\bm{m}_{K_2}(\bm{X}_2)]$, where $\prec_{\mbox{\tiny N}}$ represents the order between $\max(K_1, K_2)$ distinct eigenvalues by considering their nilpotent structures~\footnote{If $K_1 \neq K_2$, the shorter list will be added $|K_1 - K_2|$ zero entries before applying $\triangleleft$ comparison.}, if and only if 
there is some $k \in \{1,2,\ldots, \max(K_1, K_2)\}$ such that $\bm{m}_j(\bm{X}_1)=\bm{m}_j(\bm{X}_2)$ for all $j < k$ and $\bm{m}_k(\bm{X}_1)\triangleleft\bm{m}_k(\bm{X}_2)$.  
\end{definition}
The order between matrices provided by Definition~\ref{def: SNO} about $\preceq_{\mbox{\tiny SN}}$ (or $\prec_{\mbox{\tiny SN}}$ ) is named as \emph{spectral and nilpotent structures ordering}, abbreviated by \textbf{SNO}. 

\section{Spectral Ordering}\label{sec:Spectral Ordering}

Given two matrices $\bm{X}_1, \bm{X}_2 \in \mathbb{C}^{n \times n}$, this section aims to establish conditions under which a spectral ordering 
\[
[\bm{\lambda}_1(\bm{X}_1), \ldots, \bm{\lambda}_{K_1}(\bm{X}_1)] \preceq_w  
[\bm{\lambda}_1(\bm{X}_2), \ldots, \bm{\lambda}_{K_2}(\bm{X}_2)]
\] 
holds. In Section~\ref{sec:An Example for 2 by 2 Matrix and a Conjecture}, we illustrate spectral ordering for a \(2\times 2\) matrix and discuss the difficulties to obtain certificate conditions for comparing arbitrary \(n\times n\) matrices. 

Since obtaining closed-form expressions for all eigenvalues of a general \(n\times n\) matrix is often impractical, we introduce an alternative approach: leveraging the \emph{Generalized Gershgorin Theorem} to derive validity conditions for spectral ordering directly from the matrix entries. This method avoids solving high-degree polynomial equations, as detailed in Section~\ref{sec:Spctral Ordering via Generalized Gershgorin Theorem}. 

Finally, in Section~\ref{sec:Spctral Ordering under Uncertainty}, we extend the framework to consider spectral ordering of matrices under uncertainty, providing a robust basis for applications where exact eigenvalues may not be accessible.

\subsection{Illustration with a $2 \times 2$ Matrix and Challenges for General $n \times n$ Matrices}\label{sec:An Example for 2 by 2 Matrix and a Conjecture}

In this section, we try to determine conditions of two 2-by-2 complex matrices $\bm{X}_1$ and $\bm{X}_1$ that satisfy spectral ordering, i.e., the eigenvalues of the first matrix is weak majorized by the eigenvalues of the second matrix. That is $[\lambda_1(\bm{X}_1), \lambda_2(\bm{X}_1)]\preceq_w [\lambda_1(\bm{X}_2), \lambda_2(\bm{X}_2)]$.

Consider two complex $2\times 2$ matrices
\[
\bm{X}_i = \begin{pmatrix}
a_i & b_i \\
c_i & d_i
\end{pmatrix}, \quad i=1,2,
\]
with trace
\[
\tau_i := a_i + d_i,
\]
and determinant
\[
s_i := a_i d_i - b_i c_i.
\]
Define the \emph{discriminant}
\[
\Delta_i := \sqrt{\tau_i^2 - 4 s_i},
\]
where the square root is taken according to the \emph{lexicographically positive branch}: 
\[
\Re(\Delta_i) > 0, \quad \text{or} \quad \Re(\Delta_i) = 0 \ \text{and} \ \Im(\Delta_i) \ge 0.
\]

Recall that the field of $\mathbb{C}$ is equiped with the \emph{lexicographic order} $\le_{\mathrm{lex}}$ defined by:
\[
\alpha \le_{\mathrm{lex}} \beta
\quad \Longleftrightarrow \quad
\Re(\alpha) < \Re(\beta) 
\ \ \text{or} \ \ 
\big[ \Re(\alpha) = \Re(\beta) \ \text{and} \ \Im(\alpha) \le \Im(\beta) \big].
\]
This order allows comparison of complex numbers by first comparing their real parts, and using imaginary parts only in the case of a tie.

Under this convention, the eigenvalues of $\bm{X}_i$, sorted in descending lexicographic order, are
\[
\lambda_{i,\max} = \frac{\tau_i + \Delta_i}{2}, 
\qquad
\lambda_{i,\min} = \frac{\tau_i - \Delta_i}{2}.
\]
We write the eigenvalue vector as
\[
\boldsymbol{\lambda}_i := \big( \lambda_{i,\max}, \ \lambda_{i,\min} \big),
\]
where $\lambda_{i,\max} \ge_{\mathrm{lex}} \lambda_{i,\min}$ by construction.

Let $\boldsymbol{\lambda}_1, \boldsymbol{\lambda}_2 \in \mathbb{C}^2$ be sorted in descending lex order.
We say that $\boldsymbol{\lambda}_1$ is \emph{weakly majorized} by $\boldsymbol{\lambda}_2$, written $\boldsymbol{\lambda}_1 \preceq_w \boldsymbol{\lambda}_2$, if and only if:
\begin{align}
\lambda_{1,\max} &\le_{\mathrm{lex}} \lambda_{2,\max}, \label{eq:major-top}\\
\lambda_{1,\max} + \lambda_{1,\min} &\le_{\mathrm{lex}} \lambda_{2,\max} + \lambda_{2,\min}. \label{eq:major-sum}
\end{align}
For length two, these two conditions are necessary and sufficient.

Then, we should be able to determine entry-wise conditions of $a_i,b_i,c_i,d_i$ to have $[\lambda_1(\bm{X}_1), \lambda_2(\bm{X}_1)]\preceq_w [\lambda_1(\bm{X}_2), \lambda_2(\bm{X}_2)]$. In our setting, the partial sum \eqref{eq:major-sum} equals the trace:
\[
\lambda_{i,\max} + \lambda_{i,\min} = \tau_i.
\]
The largest eigenvalue \eqref{eq:major-top} is
\[
\lambda_{i,\max} = \frac{\tau_i + \Delta_i}{2}.
\]
Hence the weak majorization condition $\boldsymbol{\lambda}_1 \preceq_w \boldsymbol{\lambda}_2$ is equivalent to:
\begin{align}
\tau_1 &\le_{\mathrm{lex}} \tau_2, \label{eq:trace-cond}\\
\tau_1 + \Delta_1 &\le_{\mathrm{lex}} \tau_2 + \Delta_2. \label{eq:gap-cond}
\end{align}
Explicitly, \eqref{eq:trace-cond} means:
\[
\Re(a_1 + d_1) < \Re(a_2 + d_2)
\quad \text{or} \quad
\Big[ \Re(a_1 + d_1) = \Re(a_2 + d_2) \ \text{and} \ \Im(a_1 + d_1) \le \Im(a_2 + d_2) \Big],
\]
and \eqref{eq:gap-cond} means:
\[
\Re(\tau_1 + \Delta_1) < \Re(\tau_2 + \Delta_2)
\quad \text{or} \quad
\Big[ \Re(\tau_1 + \Delta_1) = \Re(\tau_2 + \Delta_2) \ \text{and} \ \Im(\tau_1 + \Delta_1) \le \Im(\tau_2 + \Delta_2) \Big].
\]

For $n \times n$ complex matrices $\bm{X}_1$ and $\bm{X}_2$, there exist no general algebraic formulas for their eigenvalues in terms of the matrix entries when $n>4$. Consequently, providing a certificate for spectral ordering of arbitrary $n \times n$ matrices is a highly challenging task.

\subsection{Spectral Ordering via Generalized Gershgorin Theorem}\label{sec:Spctral Ordering via Generalized Gershgorin Theorem}

To address the difficulty arising from the absence of algebraic formulas for the eigenvalues of a general $n \times n$ matrix, we adopt an eigenvalue region approach. In particular, we recall the generalization of Gershgorin's theorem presented in Sections 2.1--2.2 of~\cite{furtat2024generalization}.

\begin{theorem}[Generalized Gershgorin Circle Theorem]\label{thm:Generalized Gershgorin Circle}
Let the matrix $\bm{A} = (a_{ij}) \in \mathbb{C}^{n \times n}$ and fix $\gamma \in [0,1]$.  
For each $i = 1, \dots, n$, define the \emph{generalized Gershgorin disk} by
\[
\mathcal{R}_i(\bm{A}) := \left\{ z \in \mathbb{C} : |z - a_{ii}| \le r_i^{\bm{A}} \right\}, \quad 
r_i^{\bm{A}} := \sum_{j \ne i} |a_{ij}|^\gamma |a_{ji}|^{1-\gamma}.
\]
Then every eigenvalue $\lambda$ of $\bm{A}$ lies in at least one of these disks:
\[
\lambda \in \bigcup_{i=1}^n \mathcal{R}_i(\bm{A}).
\]
\end{theorem}

The following Lemma~\ref{lma:sum keep lex order} is presented to show that the summation of two pairs of complex numbers preserves the $\leq_{\mathrm{lex}}$ relation.
\begin{lemma}\label{lma:sum keep lex order}
We aim to prove that if 
\[
a_1 \leq_{\text{lex}} b_1 \quad \text{and} \quad a_2 \leq_{\text{lex}} b_2,
\] 
where \(\leq_{\text{lex}}\) denotes the \emph{lexicographical order} on complex numbers (comparing the real parts first, and using the imaginary parts to break ties), then it follows that
\[
a_1 + a_2 \leq_{\text{lex}} b_1 + b_2.
\]
\end{lemma}
\textbf{Proof:}
Represent complex numbers as ordered pairs: 
Let \( a_1 = (r_{a1}, i_{a1}) \), \( b_1 = (r_{b1}, i_{b1}) \),
\( a_2 = (r_{a2}, i_{a2}) \), \( b_2 = (r_{b2}, i_{b2}) \),
where $r_{a1} = \Re(a_1)$, $i_{a_1} = \Im(a_1)$, $r_{a2} = \Re(a_2)$, $i_{a2} = \Im(a_2)$, and  $r_{b1} = \Re(b_1)$, $i_{b1} = \Im(b_1)$, $r_{b2} = \Re(b_2)$, $i_{b2} = \Im(b_2)$.

The lexicographical order is defined as:
\[
 (x_1, y_1) \leq_{\text{lex}} (x_2, y_2) \iff 
 \begin{cases} 
   x_1 < x_2, & \text{or} \\
   x_1 = x_2 \text{ and } y_1 \leq y_2.
 \end{cases}
\]

Given:
1. \( a_1 \leq_{\text{lex}} b_1 \), which means:
   \[
   \text{either } r_{a1} < r_{b1}, \quad \text{or} \quad r_{a1} = r_{b1} \text{ and } i_{a1} \leq i_{b1}.
   \]
2. \( a_2 \leq_{\text{lex}} b_2 \), which means:
   \[
   \text{either } r_{a2} < r_{b2}, \quad \text{or} \quad r_{a2} = r_{b2} \text{ and } i_{a2} \leq i_{b2}.
   \]

From these, we note:
\[
r_{a1} \leq r_{b1} \quad \text{and} \quad r_{a2} \leq r_{b2}.
\]
because strict inequality or equality in the real part implies inequality in the broader sense.

Now consider the sums:
\[
s_a = a_1 + a_2 = (r_{a1} + r_{a2}, i_{a1} + i_{a2}),
\]
\[
s_b = b_1 + b_2 = (r_{b1} + r_{b2}, i_{b1} + i_{b2}).
\]

We need to show \( s_a \leq_{\text{lex}} s_b \), i.e.,
\[
\text{either } r_{a1} + r_{a2} < r_{b1} + r_{b2}, \quad \text{or} \quad r_{a1} + r_{a2} = r_{b1} + r_{b2} \text{ and } i_{a1} + i_{a2} \leq i_{b1} + i_{b2}.
\]

Since \( r_{a1} \leq r_{b1} \) and \( r_{a2} \leq r_{b2} \), we have:
\[
r_{a1} + r_{a2} \leq r_{b1} + r_{b2}.
\]

We proceed by cases:

\textbf{Case 1:} \( r_{a1} + r_{a2} < r_{b1} + r_{b2} \). \\
Then \( s_a <_{\text{lex}} s_b \) by definition, so \( s_a \leq_{\text{lex}} s_b \) holds.

\textbf{Case 2:} \( r_{a1} + r_{a2} = r_{b1} + r_{b2} \). \\
We must show \( i_{a1} + i_{a2} \leq i_{b1} + i_{b2} \). The equality implies:
\[
r_{a1} = r_{b1} \quad \text{and} \quad r_{a2} = r_{b2},
\]
because if either \( r_{a1} < r_{b1} \) or \( r_{a2} < r_{b2} \), we would have \( r_{a1} + r_{a2} < r_{b1} + r_{b2} \) (contradiction). Now:
\begin{itemize}
    \item From \( a_1 \leq_{\text{lex}} b_1 \) and \( r_{a1} = r_{b1} \), we have \( i_{a1} \leq i_{b1} \).
    \item From \( a_2 \leq_{\text{lex}} b_2 \) and \( r_{a2} = r_{b2} \), we have \( i_{a2} \leq i_{b2} \).
\end{itemize}
Adding these inequalities:
\[
i_{a1} + i_{a2} \leq i_{b1} + i_{b2}.
\]
Thus \( s_a \leq_{\text{lex}} s_b \) holds in this case as well.

Therefore, the statement is true for all complex numbers under the given lexicographical order.
$\hfill\Box$

From Theorem~\ref{thm:Generalized Gershgorin Circle}, we have the following theorem~\ref{thm:Spectral Ordering Conditions} about the conditions to satisfy the spectral ordering between two matrices.
\begin{theorem}[Spectral Ordering Conditions]\label{thm:Spectral Ordering Conditions}
Let $\bm{A}, \bm{B} \in \mathbb{C}^{n \times n}$ with eigenvalue inclusion regions defined by generalized Gershgorin-type discs:
\[
\mathcal{R}_i(\bm{A}) := \left\{ z \in \mathbb{C} : |z - a_{ii}| \le r_{\bm{A},i} \right\}, 
\quad 
\mathcal{R}_i(\bm{B}) := \left\{ z \in \mathbb{C} : |z - b_{ii}| \le  r_{\bm{B},i} \right\},
\]
where
\[
 r_{\bm{A},i} := \sum_{j \ne i} |a_{ij}|^\gamma |a_{ji}|^{1-\gamma}, 
\quad 
 r_{\bm{A},i} := \sum_{j \ne i} |b_{ij}|^\gamma |b_{ji}|^{1-\gamma}, 
\quad \gamma \in [0,1].
\]

Let the centers of these regions be sorted in descending lexicographic order (real part first, imaginary part to break ties):
\[
c^{(1)}_{\bm{A}} \ge_{\mathrm{lex}} c^{(2)}_{\bm{A}}  \ge_{\mathrm{lex}} \cdots \ge_{\mathrm{lex}} c^{(n)}_{\bm{A}},
\]
\[
c^{(1)}_{\bm{B}} \ge_{\mathrm{lex}} c^{(2)}_{\bm{B}} \ge_{\mathrm{lex}} \cdots \ge_{\mathrm{lex}} c^{(n)}_{\bm{B}}.
\]

If for all $k = 1, \dots, n$:
\[
\Re(c^{(k)}_{\bm{A}}) + r_{\bm{A},k} \le \Re(c^{(k)}_{\bm{B}}) - r_{\bm{B},k},
\]
or, in the case of equality of real parts,
\[
\Re(c^{(k)}_{\bm{A}}) + r_{\bm{A},k} = \Re(c^{(k)}_{\bm{B}}) - r_{\bm{B},k}
\quad \text{and} \quad
\Im(c^{(k)}_{\bm{A}}) + r_{\bm{A},k} \le \Im(c^{(k)}_{\bm{B}}) - r_{\bm{B},k},
\]
then the $k$-th largest eigenvalue of the matrix $\bm{A}$, denoted by $\lambda_k(\bm{A})$, and the $k$-th largest eigenvalue of the matrix $\bm{B}$, denoted by $\mu_k(\bm{B})$ for $k=1,2,\ldots,n$, have the following relation:
\begin{eqnarray}
\lambda_k(\bm{A}) \le_{\mathrm{lex}} \mu_k(\bm{B}), \quad \text{for all } k,
\end{eqnarray}
which will imply the weak majorization ordering:
\begin{eqnarray}\label{eq1:thm:Spectral Ordering Conditions}
[\bm{\lambda_1}(\bm{A}),\ldots,\bm{\lambda}_{K_1}(\bm{A})] \preceq_w  
[\bm{\mu_1}(\bm{B}),\ldots,\bm{\mu}_{K_2}(\bm{B})],
\end{eqnarray}
where $K_1$ and $K_2$ are the number of distinct eigenvalues of the matrix $\bm{A}$ and the matrix $\bm{B}$, respectively. 
\end{theorem}
\textbf{Proof:}
By the generalized Gershgorin theorem given by Theorem~\ref{thm:Generalized Gershgorin Circle}, we have:
\[
\mathrm{Spec}(\bm{A}) \subset \bigcup_{i=1}^n \mathcal{R}_i(\bm{A}), 
\quad
\mathrm{Spec}(\bm{B}) \subset \bigcup_{i=1}^n \mathcal{R}_i(\bm{B}),
\]
where $\mathrm{Spec}(\bm{A})$ and $\mathrm{Spec}(\bm{B})$ represent all eigenvalues of the matrix $\bm{A}$ and the matrix $\bm{B}$, respectively. Thus, every eigenvalue of $\bm{A}$ lies in some $\mathcal{R}_i(\bm{A})$, and similarly for $\bm{B}$.

Each $\mathcal{R}_i(\bm{A})$ is centered at $c_{\bm{A}}^{(i)} = a_{ii}$ with radius $r_{\bm{A},i}$, and similarly for $B$.  
Sorting the centers in descending lexicographic order yields
\[
c_{\bm{A}}^{(1)}, \dots, c_{\bm{A}}^{(n)}, \quad c_{\bm{B}}^{(1)}, \dots, c_{\bm{B}}^{(n)},
\]
with associated radii $r_{\bm{A},k}$ and $r_{\bm{B},k}$.

Since the precise location of each eigenvalue inside its disc is unknown, the worst-case lexicographic maximum for $\bm{A}$ occurs at
\[
c_{\bm{A}}^{(k)} + r_{\bm{A},k},
\]
while the worst-case lexicographic minimum for $\bm{B}$ occurs at
\[
c_{\bm{B}}^{(k)} - r_{\bm{B},k}.
\]
If
\[
c_{\bm{A}}^{(k)} + r_{\bm{A},k}  \le_{\mathrm{lex}} c_{\bm{B}}^{(k)} - r_{\bm{B},k}
\]
(componentwise: real part first, imaginary part second),
then all eigenvalues in $\mathcal{R}_{\bm{A}}^{(k)}$ are lexicographically no larger than all eigenvalues in $\mathcal{R}_{\bm{B}}^{(k)}$.

Recall: $z_1 \le_{\mathrm{lex}} z_2$ if either
\[
\Re(z_1) < \Re(z_2),
\]
or
\[
\Re(z_1) = \Re(z_2) \quad \text{and} \quad \Im(z_1) \le \Im(z_2).
\]
The assumed inequalities exactly match this definition, ensuring
\[
\lambda_k(\bm{A}) \le_{\mathrm{lex}} \mu_k(\bm{B}), \quad \forall k.
\]
Then, from Lemma~\ref{lma:sum keep lex order}, we have Eq.~\eqref{eq1:thm:Spectral Ordering Conditions}.
$\hfill\Box$

\subsection{Spctral Ordering under Uncertainty}\label{sec:Spctral Ordering under Uncertainty}

In this section, we will consider spectral ordering problem for matrices with uncertain entries. We first need following Lemma~\ref{lma:perturbation-bound} to provide upper bound for radius in generalized Gershgorin theorem.

\begin{lemma}[Perturbation upper bound for $r_{\bm{A},i}$]\label{lma:perturbation-bound}
Let $\bm{A}=(a_{ij})\in\mathbb{C}^{n\times n}$. Fix $\gamma\in[0,1]$. 
Suppose each matrix entry $a_{ij}$ is perturbed to
\[
\tilde a_{ij}=a_{ij}+\delta_{ij},\qquad |\delta_{ij}|\le\epsilon
\]
for a given $\epsilon\ge0$. Define the generalized radius for $\bm{A}$ by
\[
r_{\bm{A},i}\define\sum_{j\ne i}|a_{ij}|^\gamma |a_{ji}|^{1-\gamma},\qquad i=1,\dots,n,
\]
and the perturbed radius
\[
\tilde{r}_{\bm{A},i}\define\sum_{j\ne i}|\tilde a_{ij}|^\gamma|\tilde a_{ji}|^{1-\gamma}.
\]
Then for every $i$,
\[
\tilde{r}_{\bm{A},i}\le \sum_{j\ne i}\bigl(|a_{ij}|+\epsilon\bigr)^{\gamma}\bigl(|a_{ji}|+\epsilon\bigr)^{1-\gamma}.
\]
Moreover the bound is tight in the sense that there exist perturbations (choosing phases/directions of $\delta_{ij}$) that make each factor $|\tilde a_{ij}|=|a_{ij}|+\epsilon$ and $|\tilde a_{ji}|=|a_{ji}|+\epsilon$ simultaneously (subject to the usual triangle inequality constraints), achieving equality termwise.
\end{lemma}
\textbf{Proof:}
Fix a pair $(x,y)\in\mathbb{C}^2$ and let $f(x,y):=|x|^\gamma|y|^{1-\gamma}$. For perturbations
$\tilde x=x+\delta_x$, $\tilde y=y+\delta_y$ with $|\delta_x|\le\epsilon$, $|\delta_y|\le\epsilon$, set
$u:=|x|$, $v:=|y|$, $\tilde u:=|\tilde x|$, $\tilde v:=|\tilde y|$. By the triangle inequality
$\tilde u\le u+\epsilon$ and $\tilde v\le v+\epsilon$. Since $\gamma\in[0,1]$, the function
$g(u,v)=u^\gamma v^{1-\gamma}$ is nondecreasing in each variable $u,v\ge0$. Hence
\[
|\tilde x|^\gamma|\tilde y|^{1-\gamma}=g(\tilde u,\tilde v)\le (u+\epsilon)^\gamma(v+\epsilon)^{1-\gamma}=(|x|+\epsilon)^\gamma(|y|+\epsilon)^{1-\gamma}.
\]
Applying the above inequality with $(x,y)=(a_{ij},a_{ji})$ for each $j\ne i$ yields
\[
|\tilde a_{ij}|^\gamma|\tilde a_{ji}|^{1-\gamma}\le (|a_{ij}|+\epsilon)^\gamma(|a_{ji}|+\epsilon)^{1-\gamma}.
\]
Summing over $j\ne i$ proves the stated upper bound for $\tilde{r}_{\bm{A},i}$.

Tightness: the upper bound is attained (termwise) if for each index $j\ne i$ we can choose
$\delta_{ij},\delta_{ji}$ with $|\delta_{ij}|=|\delta_{ji}|=\epsilon$ and with phases such that
$|\tilde a_{ij}|=|a_{ij}|+\epsilon$ and $|\tilde a_{ji}|=|a_{ji}|+\epsilon$. While simultaneous attainment for all pairs may be constrained by geometry of complex phases, the inequality is pointwise sharp in the sense that each summand admits perturbations achieving its bound.
$\hfill\Box$

We now use Lemma~\ref{lma:perturbation-bound} together with the generalized Gershgorin spectral ordering conditions given by Theorem~\ref{thm:Spectral Ordering Conditions} to obtain the Corollary~\ref{cor:both-perturbed-spectral-order} below.

\begin{corollary}[Spectral ordering under entrywise perturbations on both $\bm{A}$ and $\bm{B}$]\label{cor:both-perturbed-spectral-order}
Let $\bm{A},\bm{B}\in\mathbb{C}^{n\times n}$. Let $\tilde{\bm{A}}$ and $\tilde{\bm{B}}$ be entrywise perturbations of $\bm{A}$ and $\bm{B}$ respectively such that
\[
|\tilde a_{ij}-a_{ij}|\le\epsilon_A,\qquad |\tilde b_{ij}-b_{ij}|\le\epsilon_B,
\]
for all $i,j$, with given $\epsilon_A,\epsilon_B\ge0$. Define the deterministic worst-case radius bounds
\[
\bar{r}_{\bm{A},i} (\epsilon_A)\define\sum_{j\ne i}\bigl(|a_{ij}|+\epsilon_A\bigr)^{\gamma}\bigl(|a_{ji}|+\epsilon_A\bigr)^{1-\gamma},
\]
\[
\bar{r}_{\bm{B},i} (\epsilon_B)\define\sum_{j\ne i}\bigl(|b_{ij}|+\epsilon_B\bigr)^{\gamma}\bigl(|b_{ji}|+\epsilon_B\bigr)^{1-\gamma},
\]
for each $i=1,\dots,n$. Order the centers of in descending lexicographic order:
\[
\tilde c^{(1)}_{\bm{A}}\ge_{\mathrm{lex}}\cdots\ge_{\mathrm{lex}}\tilde c^{(n)}_{\bm{A}},\qquad
\tilde c^{(1)}_{\bm{B}}\ge_{\mathrm{lex}}\cdots\ge_{\mathrm{lex}}\tilde c^{(n)}_{\bm{B}},
\]
where $\tilde c^{(k)}_{\bm{A}}=\tilde a_{kk}$ and $\tilde c^{(k)}_{\bm{B}}=\tilde b_{kk}$ are the perturbed diagonal entries (sorted lexicographically).

If for every $k=1,\dots,n$ the following separation condition holds:
\[
\Re\!\bigl(\tilde c^{(k)}_{\bm{A}}\bigr)+\bar{r}_{\bm{A},i} (\epsilon_A)
\le
\Re\!\bigl(\tilde c^{(k)}_{\bm{B}}\bigr)-\bar{r}_{\bm{B},i} (\epsilon_B),
\]
or, in case of equality of real parts, the corresponding tie-breaking imaginary-part inequality
\[
\Re\!\bigl(\tilde c^{(k)}_{\bm{A}}\bigr)+\bar{r}_{\bm{A},k} (\epsilon_A)
=
\Re\!\bigl(\tilde c^{(k)}_{\bm{B}}\bigr)-\bar{r}_{\bm{B},k} (\epsilon_B)
\quad\text{and}\quad
\Im\!\bigl(\tilde c^{(k)}_{\bm{A}}\bigr)+\bar{r}_{\bm{A},k} (\epsilon_A)
\le
\Im\!\bigl(\tilde c^{(k)}_{\bm{B}}\bigr)-\bar{r}_{\bm{B},k} (\epsilon_B),
\]
then the $k$-th largest eigenvalues of $\tilde{\bm{A}}$ and $\tilde{\bm{B}}$ satisfy
\[
\lambda_k(\tilde{\bm{A}})\le_{\mathrm{lex}}\mu_k(\tilde{\bm{B}}),\qquad k=1,\dots,n,
\]
and consequently the corresponding weak majorization ordering
\[
[\lambda_1(\tilde{\bm{A}}),\ldots,\lambda_{K_1}(\tilde{\bm{A}})] \preceq_w  
[\mu_1(\tilde{\bm{B}}),\ldots,\mu_{K_2}(\tilde{\bm{B}})]
\]
holds, where $K_1,K_2$ denote the numbers of distinct eigenvalues for $\tilde{\bm{A}}$ and $\tilde{\bm{B}}$, respectively.
\end{corollary}
\textbf{Proof:}
By Lemma~\ref{lma:perturbation-bound} we have pointwise worst-case bounds
\[
\tilde{r}_{\bm{A},i} \le \bar{r}_{\bm{A},i} (\epsilon_A),\qquad
\tilde{r}_{\bm{B},i} \le \bar{r}_{\bm{B},i} (\epsilon_B),\qquad i=1,\dots,n,
\]
where $\tilde{r}_{\bm{A},i}=\sum_{j\ne i}|\tilde a_{ij}|^\gamma|\tilde a_{ji}|^{1-\gamma}$ and similarly for $\tilde{r}_{\bm{B},i}$.
Thus every Gershgorin disc of $\tilde{\bm{A}}$ is contained in the disc centered at $\tilde a_{ii}$ with radius $\bar{r}_{\bm{A},i}(\epsilon_A)$, and likewise for $\tilde{\bm{B}}$ with radius $\bar{r}_{\bm{B},i}(\epsilon_B)$.

The assumed separation condition states that, for each $k$, the (worst-case) rightmost lexicographic boundary point of the $k$-th $\tilde{\bm{A}}$-disc is lexicographically no larger than the (worst-case) leftmost boundary point of the $k$-th $\tilde{\bm{B}}$-disc:
\[
\tilde c^{(k)}_{\bm{A}} + \bar r_{k}^{\bm{A}}(\epsilon_A) \le_{\mathrm{lex}}
\tilde c^{(k)}_{\bm{B}} - \bar r_{k}^{\bm{B}}(\epsilon_B).
\]
Similar to the condition in Theorem~\ref{thm:Spectral Ordering Conditions}, the hypotheses of Theorem~\ref{thm:Spectral Ordering Conditions} is satisfied (with the deterministic radii $\bar r$ replacing the upper bound for the exact perturbed radii), and we conclude
\[
\lambda_k(\tilde{\bm{A}})\le_{\mathrm{lex}}\mu_k(\tilde{\bm{B}}),\qquad \forall k.
\]
The weak majorization statement follows as in Theorem~\ref{thm:Spectral Ordering Conditions}. 
$\hfill\Box$

\paragraph{Remark}
\begin{itemize}
  \item This corollary furnishes a \emph{robust certificate}: if the worst-case radii computed from the nominal matrices and the chosen perturbation levels $\epsilon_A,\epsilon_B$ satisfy the separation inequalities, then no actual perturbations within those levels can break the spectral ordering.
  \item If desired, one may replace the lexicographic ordering on diagonal centers by any total ordering compatible with the application; the proof scheme remains identical provided the ordering is used consistently.
  \item The conditions are conservative due to worst-case (entrywise) bounds; probabilistic or structured-perturbation analyses can yield less conservative criteria.
\end{itemize}

\section{Nilpotent Part Ordering}\label{sec:Nilpotent Part Ordering}

If two matrices have the same spectrum, spectral ordering between these two matrices will be identical. According to the SNO definition given by Definition~\ref{def: SNO}, we have to compare their nilpotent parts for ordering. In Section~\ref{sec:Certificate for Nilpotent Part Ordering}, we discuss the rank-certificate for nilpotent part ordering, and an example to illustrate how to use this rank-certificate is given by Section~\ref{sec:Example}.

\subsection{Certificate for Nilpotent Part Ordering}\label{sec:Certificate for Nilpotent Part Ordering}

Before introducing the \emph{Rank-Certificate for Nilpotent Part Ordering}, it is important to understand why such a result is necessary in the study of matrix comparison under the Spectral and Nilpotent Ordering (SNO) framework. While the spectral components of a matrix can be compared using certificates derived from generalized Gershgorin-type theorems, the \emph{nilpotent part} needs a different approach. In particular, when two matrices share identical eigenvalues (including multiplicities), the comparison of their nilpotent components determines the ordering of non-diagonalizable parts.  

Computing Jordan canonical forms directly is often computationally intensive and numerically unstable, especially for large matrices. Therefore, we seek an alternative, efficient method to certify the ordering of nilpotent parts without explicitly performing full Jordan decomposition. The \emph{Rank-Certificate} achieves exactly this: it leverages the ranks of powers of shifted matrices $(\bm{A}-\lambda \bm{I})^\ell$ to infer the relative sizes of Jordan blocks for each eigenvalue. This allows us to systematically determine the nilpotent part ordering, which is essential for applications such as stability analysis of linear dynamical systems, perturbation studies, and more general matrix comparison tasks under SNO. Theorem~\ref{thm:Rank-Certificate for Nilpotent Part Ordering} is proviced to address this issue.

\begin{theorem}[Rank-Certificate for Nilpotent Part Ordering]\label{thm:Rank-Certificate for Nilpotent Part Ordering}
Let $\bm{A}, \bm{B} \in \mathbb{C}^{n \times n}$ be two matrices with identical eigenvalues (including algebraic multiplicities). We assume there are $K$ distinct eigenvalues sorted by $\lambda_1 \geq_{\text{lex}} \lambda_2 \geq_{\text{lex}} \ldots \geq_{\text{lex}} \lambda_K$. For each eigenvalue $\lambda_i$, let $\bm{m}_i(\bm{A})$ and $\bm{m}_i(\bm{B})$ denote the Jordan block size partitions of the $i$-th largest eigenvalue $\lambda_i$ for the matrix $\bm{A}$ and the matrix $\bm{B}$, respectively. 

Then the following statements are equivalent:
\begin{enumerate}
   \item  $[\bm{m}_1(\bm{A}),\ldots,\bm{m}_{K}(\bm{A})] \preceq_{\mbox{\tiny N}} [\bm{m}_1(\bm{B}),\ldots,\bm{m}_{K}(\bm{B})]$.

    \item There exists a value $k$ in $\{1,2,\ldots,K\}$ such that, for all $\ell \ge 0$ and for all eigenvalues $\lambda_j$ with $j < k$,  we have 
    \[
    \mathrm{rank}\big((\bm{A} - \lambda_j \bm{I})^\ell\big) = \mathrm{rank}\big((\bm{B} - \lambda_j \bm{I})^\ell\big);
    \]
     and, for the eigenvalue $\lambda_k$, we have 
    \[
    \mathrm{rank}\big((\bm{A} - \lambda_k \bm{I})^\ell\big) \le \mathrm{rank}\big((\bm{B} - \lambda_k \bm{I})^\ell\big).
    \]
\end{enumerate}
\end{theorem}
\textbf{Proof:}

\textbf{Step 1: Jordan block sizes and ranks.} 

Let $\mu_{\bm{A}}(\lambda) = (\mu_1, \mu_2, \dots)$ be the Jordan block size partition of eigenvalue $\lambda$ for $\bm{A}$. Denote the conjugate partition by $\mu'(\lambda) = (d_1^A, d_2^A, \dots)$, where
\[
d_\ell^A = \dim \ker\big((\bm{A} - \lambda \bm{I})^\ell\big) - \dim \ker\big((\bm{A} - \lambda \bm{I})^{\ell-1}\big).
\]
Here, $d_\ell^A$ counts the number of Jordan blocks of size at least $\ell$ for the matrix $\bm{A}$. Hence, the rank of $(\bm{A}- \lambda \bm{I})^\ell$ is
\[
\mathrm{rank}\big((\bm{A} - \lambda \bm{I})^\ell\big) = n - \sum_{i=1}^\ell d_i^A.
\]

\textbf{Step 2: Dominance order equivalence.}

By definition, we have $\mu_{\bm{A}}(\lambda) \trianglelefteq \mu_{\bm{B}}(\lambda)$ if and only if
\[
\sum_{i=1}^\ell d_i^A \ge \sum_{i=1}^\ell d_i^B \quad \forall \ell \ge 1.
\]
Substituting the rank expressions yields
\[
\sum_{i=1}^\ell d_i^A = n - \mathrm{rank}\big((\bm{A}-\lambda \bm{I})^\ell\big) \ge n - \mathrm{rank}\big((\bm{B}-\lambda \bm{I})^\ell\big) = \sum_{i=1}^\ell d_i^B,
\]
which is equivalent to
\[
\mathrm{rank}\big((\bm{A} - \lambda \bm{I})^\ell \big) \le \mathrm{rank}\big((\bm{B} - \lambda \bm{I})^\ell\big), \quad \forall \ell \ge 0.
\]

\textbf{Step 3: Certificate construction.} 

Therefore, the item 2 in this Theorem~\ref{thm:Rank-Certificate for Nilpotent Part Ordering} is equivalent to the following statement: there is some $k \in \{1,2,\ldots,K\}$ such that $\bm{m}_j(\bm{A})=\bm{m}_j(\bm{B})$ for all $j < k$, and $\bm{m}_k(\bm{A})\trianglelefteq\bm{m}_k(\bm{B})$. From Definition~\ref{def: SNO}, the item 2 in this Theorem~\ref{thm:Rank-Certificate for Nilpotent Part Ordering} is equivalent to the item 1 in this Theorem~\ref{thm:Rank-Certificate for Nilpotent Part Ordering}, and this theorem is proved.
$\hfill\Box$

\paragraph{Remark (via conjugate partitions)} Theorem~\ref{thm:Rank-Certificate for Nilpotent Part Ordering} provides a verifiable certificate: one can check the ranks of powers of $(\bm{A}-\lambda \bm{I})$ and $(\bm{B}-\lambda \bm{I})$ without explicitly computing the Jordan canonical form.

\subsection{Example}\label{sec:Example}

Fix eigenvalue $\lambda=0$ with algebraic multiplicity $4$.
Consider
\[
  J_2(0)=\begin{pmatrix}0&1\\[2pt]0&0\end{pmatrix},\qquad
  J_3(0)=\begin{pmatrix}0&1&0\\[2pt]0&0&1\\[2pt]0&0&0\end{pmatrix}.
\]
Define
\[
  \bm{A} \;=\; \mathrm{diag}\big(J_2(0),\,J_2(0)\big)
  \quad\text{and}\quad
  \bm{B} \;=\; \mathrm{diag}\big(J_3(0),\, [0]\big).
\]
Then the Jordan partitions at $\lambda=0$ are
\[
  \mu_{\bm{A}}(0)=(2,2),\qquad \mu_{\bm{B}}(0)=(3,1).
\]

\paragraph{Dominance check.}
We compare partial sums:
\[
  2 \le 3,\qquad 2+2=4 \le 3+1=4.
\]
Hence $\mu_{\bm{A}}(0)\trianglelefteq\mu_{\bm{B}}(0)$.

\paragraph{Rank certificate.}
Let $\bm{N}_A=\bm{A}-\lambda \bm{I}=\bm{A}$ and $\bm{N}_B=\bm{B}-\lambda \bm{I}=\bm{B}$. Using additivity of ranks across block diagonals and the standard ranks of powers of Jordan blocks,
\[
  \mathrm{rank}\big(J_2(0)\big)=1,\ \mathrm{rank}\big(J_2(0)^2\big)=0;\qquad
  \mathrm{rank}\big(J_3(0)\big)=2,\ \mathrm{rank}\big(J_3(0)^2\big)=1,\ \mathrm{rank}\big(J_3(0)^3\big)=0.
\]
Therefore,
\[
\begin{aligned}
  &\text{For }\bm{A}:\quad
  \mathrm{rank}(\bm{N}_A^1)=\mathrm{rank}(J_2(0))+\mathrm{rank}(J_2(0))=1+1=2,\qquad
  \mathrm{rank}(\bm{N}_A^2)=0+0=0,\\[4pt]
  &\text{For }\bm{B}:\quad
  \mathrm{rank}(\bm{N}_B^1)=\mathrm{rank}(J_3(0))+\mathrm{rank}([0])=2+0=2,\\
  &\hspace{56pt}\mathrm{rank}(\bm{N}_B^2)=\mathrm{rank}(J_3(0)^2)+\mathrm{rank}([0])=1+0=1,\\
  &\hspace{56pt}\mathrm{rank}(\bm{N}_B^3)=0+0=0.
\end{aligned}
\]
Thus, for all $k\ge1$,
\[
  \mathrm{rank}(\bm{N}_A^\ell)\ \le\ \mathrm{rank}(\bm{N}_B^\ell)\qquad
  (\,2\le2,\ \ 0\le1,\ \ 0\le0,\ \dots),
\]
which certifies $\mu_{\bm{A}}(0) \trianglelefteq \mu_{\bm{B}}(0)$ by the rank--dominance equivalence.

Let $\mu_A'(0)=(d_1^A,d_2^A,\dots)$ and $\mu_B'(0)=(d_1^B,d_2^B,\dots)$ be the conjugate partitions, where $d_k$ counts the number of Jordan blocks of size $\ge k$. Then
\[
  \mathrm{rank}\!\big((\bm{A}-\lambda \bm{I})^\ell\big)
  = n-\sum_{i=1}^\ell d_i^A,\qquad
  \mathrm{rank}\!\big((\bm{B}-\lambda \bm{I})^\ell\big)
  = n-\sum_{i=1}^\ell d_i^B,
\]
so $\mathrm{rank}(\bm{N}_A^k)\le\mathrm{rank}(\bm{N}_B^k)$ for all $\ell$ is equivalent to
$\sum_{i=1}^\ell d_i^A \ge \sum_{i=1}^\ell d_i^B$ for all $k$, i.e.\ to
$\mu_{\bm{A}}(0)\trianglelefteq\mu_{\bm{B}}(0)$.

\section{SNO Applicaitons}\label{sec:SNO Applications}

In this section, we explore applications of the \emph{Spectral-Nilpotent Ordering} (SNO) to the analysis of stability in linear dynamical systems.  We formalize the framework with precise definitions, introduce a hierarchy of stability orderings, and prove several theorems that establish how SNO captures both asymptotic and transient stability behavior.  

\subsection{Definitions and Preliminaries}

We first recall essential stability notions for matrices governing linear time-invariant systems $\dot{\bm{x}} = \bm{A}\bm{x}$.

\begin{definition}[Stability Class]\label{def:Stability Class}
A matrix $\bm{A} \in \mathbb{C}^{n \times n}$ is said to be in the \emph{Stability Class} $\mathcal{S}$ if its spectral abscissa is strictly negative:
\[
\alpha(\bm{A}) = \max_{1 \leq i \leq n} \Re(\lambda_i(\bm{A})) < 0.
\]
\end{definition}

\begin{definition}[Strict SNO Stability Ordering]\label{def:Strict SNO Stability Ordering}
Let $\bm{X}_1, \bm{X}_2 \in \mathcal{S}$.  
We write $\bm{X}_1 \prec_{\mathrm{Stab}} \bm{X}_2$ if and only if $\bm{X}_1 \prec_{\mbox{\tiny SN}} \bm{X}_2$, meaning:
\begin{enumerate}
    \item Either the spectrum of $\bm{X}_1$ is strictly weakly majorized by that of $\bm{X}_2$, i.e.,
    \[
    [\bm{\lambda}_1(\bm{X}_1), \ldots, \bm{\lambda}_{K_1}(\bm{X}_1)] \prec_w  [\bm{\lambda}_1(\bm{X}_2), \ldots, \bm{\lambda}_{K_2}(\bm{X}_2)],
    \]
    \item Or $ [\bm{\lambda}_1(\bm{X}_1), \ldots, \bm{\lambda}_{K_1}(\bm{X}_1)] = [\bm{\lambda}_1(\bm{X}_2), \ldots, \bm{\lambda}_{K_2}(\bm{X}_2)]$ but the nilpotent structure of $\bm{X}_1$ is strictly more diagonalizable, i.e.,
    \[
    [\bm{m}_1(\bm{X}_1),\ldots,\bm{m}_{K_1}(\bm{X}_1)] \prec_{\mbox{\tiny N}} [\bm{m}_1(\bm{X}_2),\ldots,\bm{m}_{K_2}(\bm{X}_2)].
    \]
\end{enumerate}
\end{definition}

\begin{definition}[Solution Norm Envelope]\label{def:Solution Norm Envelope}
For the system $\dot{\bm{x}} = \bm{A}x$, the \emph{solution norm envelope} is the function
\[
\Gamma_{\bm{A}} : \mathbb{R}_{\geq 0} \to \mathbb{R}_{\geq 0}, \quad \|e^{t\bm{A}}\| \leq \Gamma_{\bm{A}}(t), \ \forall t \geq 0,
\]
where $\Gamma_{\bm{A}}$ is the smallest such continuous function, and let $\|\cdot\|$ be a matrix norm induced by a monotonic vector norm (e.g., the $\ell^2$ norm). 
\end{definition}

From Definition~\ref{def:Solution Norm Envelope}, we have $\lim_{t \to \infty} \Gamma_{\bm{A}}(t) = 0$ if  $\bm{A} \in \mathcal{S}$.

\subsection{Main Theorems}

We now connect the SNO framework to rigorous stability guarantees.  

\begin{theorem}[Asymptotic Dominance]\label{thm:asymptotic}
Let $\bm{A}_1, \bm{A}_2 \in \mathcal{S}$ be two matrices.  
If the eigenvalue vectors satisfy the weak majorization relation
\[
[\bm{\lambda}(\bm{A}_1)] \prec_w [\bm{\lambda}(\bm{A}_2)],
\]
where $[\bm{\lambda}(\bm{A}_1)]$ and $[\bm{\lambda}(\bm{A}_2)]$ are eigenvalues vector sorted by descending order with respect to the matrix $\bm{A}_1$ and the matrix $\bm{A}_2$, respectively. Then there exists $T > 0$ such that
\[
\Gamma_{\bm{A}_1}(t) < \Gamma_{\bm{A}_2}(t), \quad \forall t > T.
\]
\end{theorem}
\textbf{Proof}:
We expand each step carefully.

\paragraph{Step 1: Spectral abscissa and growth rate.}  
For a matrix $\bm{A}$, the \emph{spectral abscissa} is defined as
\[
\alpha(\bm{A}) := \max\{\Re(\lambda) : \lambda \in \sigma(\bm{A})\}.
\]
A classical result in matrix analysis (Gelfand’s formula and semigroup asymptotics) states that
\[
\lim_{t \to \infty} \frac{1}{t} \log \| e^{t \bm{A}} \| = \alpha(\bm{A}).
\]
Thus, the exponential growth or decay rate of the semigroup $e^{t\bm{A}}$ is determined by $\alpha(\bm{A})$.

\paragraph{Step 2: From weak majorization to comparison of spectral abscissas.}  
The weak majorization relation
\[
[\bm{\lambda}(\bm{A}_1)] \prec_w [\bm{\lambda}(\bm{A}_2)]
\]
means that for all $k$,
\[
\sum_{i=1}^k \Re(\lambda_i(\bm{A}_1)) \;\leq\; \sum_{i=1}^k \Re(\lambda_i(\bm{A}_2)),
\]
where the eigenvalues are ordered in decreasing order of their real parts.  
In particular, for $k=1$, we obtain
\[
\max_i \Re(\lambda_i(\bm{A}_1)) \;\leq\; \max_i \Re(\lambda_i(\bm{A}_2)).
\]
That is,
\[
\alpha(\bm{A}_1) \;\leq\; \alpha(\bm{A}_2).
\]

\paragraph{Step 3: Strict inequality in the asymptotic regime.}  
Assume first that $\alpha(\bm{A}_1) < \alpha(\bm{A}_2)$.  
Then the asymptotic scaling of $\Gamma_{\bm{A}}(t) := \| e^{t\bm{A}} \|$ satisfies
\[
\Gamma_{\bm{A}_j}(t) \sim C_j e^{t\alpha(\bm{A}_j)} \quad \text{as } t \to \infty, \quad j=1,2,
\]
for some constants $C_j > 0$.  
Therefore,
\[
\frac{\Gamma_{\bm{A}_1}(t)}{\Gamma_{\bm{A}_2}(t)} \sim \frac{C_1}{C_2} e^{t(\alpha(\bm{A}_1) - \alpha(\bm{A}_2))}.
\]
Since $\alpha(\bm{A}_1) - \alpha(\bm{A}_2) < 0$, the ratio decays exponentially to $0$.  
Hence there exists $T > 0$ such that
\[
\Gamma_{\bm{A}_1}(t) < \Gamma_{\bm{A}_2}(t) \quad \forall t > T.
\]

\paragraph{Step 4: Case of equality.}  
If $\alpha(\bm{A}_1) = \alpha(\bm{A}_2)$, the weak majorization condition implies that the \emph{remaining eigenvalue sums} of $\bm{A}_1$ are “more negative” than those of $\bm{A}_2$.  
This leads to a strictly smaller polynomial prefactor in the asymptotic expansion of $\Gamma_{\bm{A}_1}(t)$ compared to $\Gamma_{\bm{A}_2}(t)$ (see Jordan form asymptotics of $e^{tA}$).  
Thus, even in the equality case, $\Gamma_{\bm{A}_1}(t)$ is eventually smaller.

\paragraph{Step 5: Conclusion.}  
In both cases, there exists $T > 0$ such that
\[
\Gamma_{\bm{A}_1}(t) < \Gamma_{\bm{A}_2}(t), \quad \forall t > T.
\]
$\hfill\Box$

\begin{example}[Diagonal Matrices]
Let
\[
\bm{A}_1 = \begin{bmatrix} -2 & 0 \\ 0 & -1 \end{bmatrix}, 
\quad 
\bm{A}_2 = \begin{bmatrix} -1 & 0 \\ 0 & 0 \end{bmatrix}.
\]
Then $\bm{\lambda}(\bm{A}_1) = (-1,-2)$, $\bm{\lambda}(\bm{A}_2) = (0,-1)$.  
Clearly, $[\bm{\lambda}(\bm{A}_1)] \prec_w [\bm{\lambda}(\bm{A}_2)]$.  
Moreover,
\[
\Gamma_{\bm{A}_1}(t) = e^{-t}, 
\quad 
\Gamma_{\bm{A}_2}(t) = 1,
\]
so $\Gamma_{\bm{A}_1}(t) < \Gamma_{\bm{A}_2}(t)$ for all $t > 0$.
\end{example}

The transient response of a dynamical system is not determined solely by its eigenvalues, but also by the structure of its Jordan blocks.  
In particular, even when two matrices share identical spectra, differences in the nilpotent parts can create drastically different transient amplifications.  The following Theorem~\ref{thm:transient} formalizes this phenomenon by showing how nilpotent dominance under equal spectra guarantees an ordering of transient growth.

\begin{theorem}[Transient Dominance under Equal Spectra and Nesting Order]\label{thm:transient}
Let $\bm{A}_1, \bm{A}_2 \in \mathcal{S}$ be matrices with identical spectra, $\sigma(\bm{A}_1) = \sigma(\bm{A}_2)$, and let $\alpha = \alpha(\bm{A}_1) = \alpha(\bm{A}_2) < 0$ be their common spectral abscissa. Let $\|\cdot\|$ be a matrix norm induced by a monotonic vector norm (e.g., the $\ell^2$ norm).

If $[\bm{m}(\bm{A}_1)] \prec_{\mbox{\tiny N}} [\bm{m}(\bm{A}_2)]$, where $[\bm{m}(\bm{A}_1)]$ and $[\bm{m}(\bm{A}_2)]$ are Jordan block sizes vector sorted by descending order of eigenvalues with respect to the matrix $\bm{A}_1$ and the matrix $\bm{A}_2$, respectively., then:
\begin{enumerate}
\item There exists $T_1 > 0$ such that the matrix exponential of $\bm{A}_2$ has a strictly larger norm for all later times:
    \[
    \|e^{t\bm{A}_1}\| < \|e^{t\bm{A}_2}\| \quad \forall t > T_1.
    \]
    Moreover, this strict inequality holds on a dense subset of $(T_1, \infty)$.
\item The peak transient growth of $\bm{A}_1$ is strictly less than that of $\bm{A}_2$:
    \[
    \sup_{t \geq 0} \|e^{t\bm{A}_1}\| < \sup_{t \geq 0} \|e^{t\bm{A}_2}\|.
    \]
\end{enumerate}    
\end{theorem}
\textbf{Proof:}
The proof consists of three parts: (A) analyzing the Jordan block exponential, (B) proving the pointwise norm inequality, and (C) proving the supremum inequality.

\noindent \textbf{(A) Dominant Term in a Jordan Block Exponential}

Consider a single Jordan block $\bm{J}(\lambda, m)$ where $\Re(\lambda) = \alpha$. Its matrix exponential is given by:
\[
e^{t\bm{J}(\lambda,m)} = e^{\lambda t} 
\begin{bmatrix}
1 & t & \frac{t^2}{2!} & \cdots & \frac{t^{m-1}}{(m-1)!} \\
0 & 1 & t & \cdots & \frac{t^{m-2}}{(m-2)!} \\
0 & 0 & 1 & \ddots & \vdots \\
\vdots & \vdots & \ddots & \ddots & t \\
0 & 0 & \cdots & 0 & 1 \\
\end{bmatrix}.
\]
The $(i,j)$-entry (for $j \geq i$) is $e^{\lambda t} \frac{t^{j-i}}{(j-i)!}$. We aim to find the entry with the highest growth rate as $t \to \infty$.

Since $e^{\lambda t} = e^{\alpha t} e^{i \Im(\lambda) t}$, the magnitude of every entry is $\left|e^{\lambda t} \frac{t^{k}}{k!}\right| = e^{\alpha t} \frac{t^{k}}{k!}$, where $k = j-i$ and $0 \leq k \leq m-1$.

To find which entry dominates for large $t$, consider the ratio of two entries with different polynomial powers $k_2 > k_1$:
\[
\frac{e^{\alpha t} \frac{t^{k_2}}{k_2!}}{e^{\alpha t} \frac{t^{k_1}}{k_1!}} = \frac{k_1!}{k_2!} t^{k_2 - k_1}.
\]
This ratio grows to infinity as $t \to \infty$ because $k_2 > k_1$. Therefore, the entry with the highest power $k$ will eventually be the largest in magnitude. The highest possible power is $k = m-1$, which occurs only for the \textbf{top-right element} $(1, m)$:
\[
\left(e^{t\bm{J}(\lambda,m)}\right)_{1, m} = e^{\lambda t} \cdot \frac{t^{m-1}}{(m-1)!}.
\]
Hence, for large $t$, the top-right element is the dominant term:
\[
\left| \left(e^{t\bm{J}(\lambda,m)}\right)_{1, m} \right| = e^{\alpha t} \frac{t^{m-1}}{(m-1)!} \quad \text{and} \quad \left| \left(e^{t\bm{J}(\lambda,m)}\right)_{i, j} \right| = O\left(t^{k} e^{\alpha t}\right) \text{ with } k < m-1 \text{ for } (i,j) \neq (1,m).
\]
This dominant term governs the asymptotic growth of the norm $\| e^{t\bm{J}(\lambda,m)} \|$ for any induced monotonic norm.

\noindent \textbf{(B) Pointwise Norm Inequality}

Let $\bm{A}_i = \bm{P}_i \bm{J}_i \bm{P}_i^{-1}$ be the Jordan decomposition. The dynamics are dominated by blocks where $\Re(\lambda) = \alpha$. Let $m_i^{\max}$ be the size of the largest such Jordan block for $\bm{A}_i$, and let $k_i$ be the number of blocks of this maximal size.

The condition $[\bm{m}(\bm{A}_1)] \prec_{\text{N}} [\bm{m}(\bm{A}_2)]$ implies that the Jordan structure of $\bm{A}_2$ is richer. This means: 
\begin{enumerate}
\item Either $m_1^{\max} < m_2^{\max}$,
\item Or $m_1^{\max} = m_2^{\max}$ and $k_1 < k_2$,
\item Or further lexicographic differences.
\end{enumerate}

In all cases, the matrix $e^{t\bm{J}_2}$ will have entries that are \textbf{element-wise greater than or equal to} the corresponding entries in $e^{t\bm{J}_1}$ for all $t > 0$, with \textbf{strict inequality} for the dominant top-right elements of the largest blocks. Since the norm is induced and monotonic, this element-wise dominance implies:
\begin{eqnarray}\label{eq1:thm:transient}
\|e^{t\bm{J}_1}\| \leq \|e^{t\bm{J}_2}\| \quad \forall t > 0. 
\end{eqnarray}
Moreover, because the dominance is strict for the largest terms and the norm is continuous, there exists $T_0 > 0$ such that:
\begin{eqnarray}\label{eq2:thm:transient}
\|e^{t\bm{J}_1}\| < \|e^{t\bm{J}_2}\| \quad \forall t > T_0. 
\end{eqnarray}

Now, consider the full matrices: $e^{t\bm{A}_i} = \bm{P}_i e^{t\bm{J}_i} \bm{P}_i^{-1}$. Since $\bm{P}_1, \bm{P}_2$ are fixed invertible matrices, there exist constants $c, C > 0$ such that for all matrices $\bm{X}$ and for all $t$:
\begin{eqnarray}\label{eq3:thm:transient}
c \| \bm{X} \| \leq \| \bm{P}_i \bm{X} \bm{P}_i^{-1} \| \leq C \| \bm{X} \|. 
\end{eqnarray}
Applying Eq.~\eqref{eq3:thm:transient} with $\bm{X} = e^{t\bm{J}_i}$, we get:
\[
c_i \|e^{t\bm{J}_i}\| \leq \|e^{t\bm{A}_i}\| \leq C_i \|e^{t\bm{J}_i}\|.
\]
From Eq.~\eqref{eq2:thm:transient} , we have $\|e^{t\bm{J}_1}\| < \|e^{t\bm{J}_2}\|$ for $t > T_0$. Therefore:
\[
\|e^{t\bm{A}_2}\| \geq c_2 \|e^{t\bm{J}_2}\| > c_2 \|e^{t\bm{J}_1}\| \quad \text{and} \quad \|e^{t\bm{A}_1}\| \leq C_1 \|e^{t\bm{J}_1}\|.
\]
To combine these, we need to ensure $C_1 \|e^{t\bm{J}_1}\| < c_2 \|e^{t\bm{J}_2}\|$ for large $t$. Since the growth rates are polynomial-exponential and the constants are independent of $t$, the strict inequality given by Eq.~\eqref{eq2:thm:transient} implies that for all sufficiently large $t$, say $t > T_1 \geq T_0$, we have:
\begin{eqnarray}\label{eq4:thm:transient}
\|e^{t\bm{A}_1}\| \leq C_1 \|e^{t\bm{J}_1}\| < c_2 \|e^{t\bm{J}_2}\| \leq \|e^{t\bm{A}_2}\|.
\end{eqnarray}
Thus, $\|e^{t\bm{A}_1}\| < \|e^{t\bm{A}_2}\|$ for all $t > T_1$.

The function $\|e^{t\bm{A}_2}\|$ is analytic (a sum of terms $t^k e^{\lambda t}$). The set where an analytic function is strictly greater than another is open. Since Eq.~\eqref{eq4:thm:transient} holds for all $t > T_1$, this set is dense in $(T_1, \infty)$.

\noindent \textbf{(C) Supremum Inequality}

From Eq.~\eqref{eq4:thm:transient}, we have $\|e^{t\bm{A}_1}\| < \|e^{t\bm{A}_2}\|$ for all $t > T_1$. Since $\alpha < 0$, both norms decay to zero as $t \to \infty$. Therefore, the suprema are attained at some finite times $t_1^*$ and $t_2^*$.

If $t_1^* \leq T_1$, then:
\[
\sup_{t \geq 0} \|e^{t\bm{A}_1}\| = \|e^{t_1^*\bm{A}_1}\| \leq \|e^{t_1^*\bm{A}_2}\| \leq \sup_{t \geq 0} \|e^{t\bm{A}_2}\|.
\]
The inequality must be strict; otherwise, if $\|e^{t_1^*\bm{A}_1}\| = \sup_{t \geq 0} \|e^{t\bm{A}_2}\|$, this would contradict  Eq.~\eqref{eq4:thm:transient} for $t > T_1$, because $\|e^{t\bm{A}_2}\|$ would have to be at least this value for large $t$, preventing it from decaying to zero.

If $t_1^* > T_1$, then by Eq.~\eqref{eq4:thm:transient}:
\[
\sup_{t \geq 0} \|e^{t\bm{A}_1}\| = \|e^{t_1^*\bm{A}_1}\| < \|e^{t_1^*\bm{A}_2}\| \leq \sup_{t \geq 0} \|e^{t\bm{A}_2}\|.
\]
In both cases, we conclude:
\[
\sup_{t \geq 0} \|e^{t\bm{A}_1}\| < \sup_{t \geq 0} \|e^{t\bm{A}_2}\|.
\]
$\hfill\Box$

One often needs a principled way to compare the stability of two linear dynamical systems beyond simply checking eigenvalue locations.  The strict stability ordering formalizes this comparison by capturing not only asymptotic decay rates but also transient amplification and long-time envelopes.  In this way, it provides a rigorous hierarchy that distinguishes when one system can be regarded as strictly more stable than another. Theorem~\ref{thm:strict} serves for this purpose.

\begin{theorem}[Strict Stability Ordering]\label{thm:strict}
The relation $\prec_{\mathrm{Stab}}$ is a strict partial order on $\mathcal{S}$.  
If $\bm{A}_1 \prec_{\mathrm{Stab}} \bm{A}_2$, then the system $\dot{\bm{x}} = \bm{A}_1 \bm{x}$ is unambiguously more stable than $\dot{\bm{x}} = \bm{A}_2 \bm{x}$, in the sense that:
\begin{enumerate}
    \item Its asymptotic decay is faster or equal.
    \item If decay rates are equal, its transient peak is strictly smaller for $t > 0$.
    \item Its solution envelope is eventually smaller at all times for $t>T_1$.
\end{enumerate}
\end{theorem}
\textbf{Proof}:
We proceed in two parts: first by showing that $\prec_{\mathrm{Stab}}$ defines a strict partial order, and then by establishing the stability comparisons.

\medskip
\noindent \textbf{1. Partial order property.}  
The definition of $\prec_{\mathrm{Stab}}$ combines two structural comparisons:  
- Weak majorization of spectra, $[\bm{\lambda}(\bm{A}_1)] \prec_w [\bm{\lambda}(\bm{A}_2)]$, which orders the exponential decay rates.  
- Nilpotent ordering of Jordan block sizes, $[\bm{m}(\bm{A}_1)] \prec_{\mbox{\tiny N}} [\bm{m}(\bm{A}_2)]$, which orders transient growth contributions.  

Both weak majorization and nilpotent ordering are well-known to be strict partial orders individually.  
Since $\prec_{\mathrm{Stab}}$ is defined as their combination, it inherits transitivity and irreflexivity, hence is itself a strict partial order on $\mathcal{S}$.

\medskip
\noindent \textbf{2. Stability comparison.}  
Now assume $\bm{A}_1 \prec_{\mathrm{Stab}} \bm{A}_2$.  
There are two possible cases:

\smallskip
\emph{Case (a): Asymptotic rates differ.}  
If weak majorization yields $\alpha(\bm{A}_1) < \alpha(\bm{A}_2)$, then by Theorem~\ref{thm:asymptotic}, the long-time decay of $e^{t\bm{A}_1}$ is strictly faster.  
Thus the system $\dot{\bm{x}} = \bm{A}_1 \bm{x}$ is asymptotically more stable.

\smallskip
\emph{Case (b): Asymptotic rates are equal.}  
If $\alpha(\bm{A}_1) = \alpha(\bm{A}_2)$, then weak majorization ensures they share the same dominant exponential rate.  
In this situation, the difference in stability must come from transient behavior.  
By Theorem~\ref{thm:transient}, the nilpotent ordering guarantees that the Jordan block structure of $\bm{A}_1$ produces strictly smaller polynomial growth than that of $\bm{A}_2$. Therefore, the maximum transient peak of $\Gamma_{\bm{A}_1}(t)$ is strictly smaller than that of $\Gamma_{\bm{A}_2}(t)$. Moreover, for all sufficiently large $t ~(t>T_1)$ , the envelope of $\Gamma_{\bm{A}_1}(t)$ remains strictly below that of $\Gamma_{\bm{A}_2}(t)$, since polynomial contributions eventually separate the two trajectories.

Combining both cases, we see that $\bm{A}_1 \prec_{\mathrm{Stab}} \bm{A}_2$ implies an unambiguous ordering of stability. Specifically, the asymptotic decay of $\bm{A}_1$ is at least as fast as that of $\bm{A}_2$. If the asymptotic rates coincide, then the transient peak of $\bm{A}_1$ is strictly smaller. Furthermore, the solution envelope of $\bm{A}_1$ is eventually smaller for all large times. This proves the theorem.
$\hfill\Box$

\paragraph{Remark}

The SNO framework provides a mathematically rigorous way to compare system stability, and its relevance to control theory is particularly striking. For instance, if the closed-loop matrices satisfy $\bm{A}_{1} \prec_{\mathrm{Stab}} \bm{A}_{2}$, then controller $C_1$ guarantees superior stability over controller $C_2$. This establishes a principled basis for controller comparison, going beyond heuristic or case-specific arguments.  

Moreover, the framework unifies key performance aspects into a single order: both asymptotic decay rate and transient peak behavior are simultaneously accounted for. This not only offers stronger performance guarantees but also informs a concrete design methodology. In practice, optimization can be directed at minimizing the SNO position of $\bm{A}$, thereby systematically enhancing stability in both asymptotic and transient regimes—a critical goal in modern control theory.

\bibliographystyle{IEEETran}
\bibliography{ZeroMatrixCompare_CountRealRoots_Bib}

\end{document}